\documentclass[11pt,a4paper,leqno]{amsart}

\usepackage{amsmath,,amssymb,amsthm}
\usepackage[english]{babel}
\usepackage{xcolor}
\usepackage{datetime}
\usepackage{enumerate}
\usepackage[foot]{amsaddr}
\usepackage{lipsum}

\usepackage[colorlinks=true,
    linkcolor=blue, citecolor=blue,urlcolor=blue, linktocpage=true]{hyperref}

\newtheorem{theorem}{Theorem}

\newtheorem{theo}{Theorem}[section]

\newtheorem{propo}[theo]{Proposition}

\newtheorem{lem}[theo]{Lemma}

\theoremstyle{definition}
\newtheorem{defin}[theo]{Definition}

\newtheorem{remark}[theo]{Remark}

\newtheorem{question}{Question}

\numberwithin{equation}{section}

\usepackage[noadjust]{marginnote}
\def\R{\mathbb{R}}
\def\C{\mathbb{C}}
\def\D{\mathbb{D}}
\def\Z{\mathbb{Z}}

\def\E{\mathbf{E}}
\def\V{\mathbf{V}}
\def\P{\mathbf{P}}

\def\K{\mathcal{K }}

\title{Khinchin families, set constructions, partitions and exponentials}

\author[A. Cant\'{o}n]{Alicia Cant\'{o}n}

\address[Alicia Cant\'{o}n]{Departamento de Matem\'{a}tica e
Inform\'{a}tica Aplicadas a las Ingenier\'{\i}as Civil y Naval, ETSIN,
Universidad Polit\'{e}cnica de Madrid, Avenida de la Memoria~4, Ciudad Universitaria, 28040 Madrid, Spain}
\email{alicia.canton@upm.es}

\author[J.\,L. Fern\'{a}ndez]{Jos\'{e} L. Fern\'{a}ndez}
\address[Jos\'{e} L. Fern\'{a}ndez, Pablo Fern\'{a}ndez and V\'{\i}ctor J. Maci\'{a}]{Departamento de Matem\'{a}ticas, Facultad de Ciencias, Universidad Aut\'{o}noma de Madrid, {Ciudad Universitaria de Cantoblanco s/n, 28049 Madrid}, Spain.}
\email{joseluis.fernandez@uam.es}

\author[P. Fern\'{a}ndez]{Pablo Fern\'{a}ndez}
\email{pablo.fernandez@uam.es}

\author[V.\,J. Maci\'{a}]{V\'{\i}ctor J. Maci\'{a}}
\email{victor.macia@uam.es}



\begin{document}

\begin{abstract}In this paper, 
 we give a simple criterion to verify that functions of the form $e^g$ are in the Hayman class when $g$ is a power series with nonnegative coefficients. Thus, using  the Hayman and B\'{a}ez-Duarte formulas,
 we obtain asymptotics for the coefficients of generating functions that arise in many examples of set construction in analytic combinatorics. This new criterion greatly simplifies that obtained previously by the authors.
\end{abstract}

\maketitle

\section{Introduction}\label{sec1}

The present paper is a follow-up to \cite{K_uno}. There,  following the lead of Hayman~\cite{Hayman}, Rosenbloom \cite{Rosenbloom} and B\'{a}ez-Duarte \cite{BaezDuarte}, a basic theory of the so called Khinchin families is laid out. The notions of Gaussian and strongly Gaussian power series (which include power series in the Hayman class) are developed, and it is shown how the asymptotic formulas of Hayman and B\'{a}ez-Duarte (see Theorem \ref{eq:condition baez-duarte}) provide a convenient way of handling the asymptotic of coefficients of strongly Gaussian power series.

\smallskip
In this context, a particularly interesting question consists of establishing the asymptotic behavior of the coefficients of a power series $f$ that is written as $f= e^g$, where $g$ is a power series with nonnegative coefficients.

For instance, in Combinatorics, the \textit{set construction}, both labeled and unlabeled,  is concerned with the combinatorial class of sets formed with  objects drawn from a given combinatorial  class. The generating function of the class of sets is of the form $f=e^g$, where \textit{typically} $g$ is a power series with nonnegative coefficients. The generating function of the Bell numbers and the partition function are examples of functions arising from the set construction. See Section \ref{seccion:set construction} of the present paper, or consult the comprehensive treatise \cite{Flajolet} of Flajolet and  Sedgewick.

Hayman, in \cite{Hayman}, deals with this question in the following particular cases: when $g$ (not necessarily with nonnegative coefficients) already belongs to the Hayman class, and when $g$ is a polynomial with nonnegative coefficients that obey certain divisibility properties. See the beginning of Section~\ref{seccion:exponentials and Hayman} for details.

In Theorem 4.1 of \cite{K_uno},  a basic criterion is presented that can be used to establish conditions on the power series $g$ with nonnegative coefficients which imply that $f=e^g$ is in the Hayman class and, therefore, is strongly Gaussian. This criterion, combined with the Hayman and B\'{a}ez-Duarte asymptotic formulas, gives asymptotic results for the coefficients of a large and varied collection of power series.

In this paper, this criterion is greatly simplified and its applicability expanded.  Theorems \ref{teor:general for quasiexponentials} and \ref{teor:general for quasigeometric}, the main results of the present paper, exhibit conditions on the (nonnegative) coefficients of a power series $g$ which guarantee that $e^g$ is in the Hayman class. The criteria mentioned above are written directly in terms of the \emph{growth} of the coefficients of $g$. Compare with Theorem 4.1 in \cite{K_uno}, that is reproduced as Theorem~\ref{teor:cut for e^g in Hayman class} in this paper.

By appealing to Theorems \ref{teor:general for quasiexponentials} and \ref{teor:general for quasigeometric}, the verification that the generating functions of many combinatorial set constructions is in the Hayman class becomes immediate; asymptotic formulas for  their coefficients follow  simply from the Hayman and B\'{a}ez-Duarte asymptotic formulas.

Although this paper is a follow-up to \cite{K_uno}, we intend this presentation to be self-contained. The most salient  and relevant features of the theory of Khinchin families are described within the present paper; nonetheless, we refer to \cite{K_uno} for a detailed treatment.

\subsection{Notation and some preliminaries}

The symbol $a_n\sim b_n$ as $n\to\infty$ means that $a_n/b_n\to 1$ as $n\to\infty$, while $a_n\asymp b_n$ as $n\to\infty$ means that $1/C\le a_n/b_n<C$ for some $C>0$.

We use $\E$, $\V$ and $\P$ to denote expectation, variance and probability generically.

For random variables $X $ and $Y$, we write $X\stackrel{d}{=}Y$ to signify  that $X$ and $Y$ have the same distribution.

If $(Z_n)_{n \ge 1}$ is a sequence  of random variables and $Z$ is another random variable, the notation
$Z_n \stackrel{\rm d}{\longrightarrow} Z$ as $n\to\infty$ means convergence in distribution, which by Levy's theorem is equivalent to pointwise convergence of characteristic functions, i.e.,
\[\lim_{n \to \infty} \E(e^{ \imath \theta Z_n})=\E(e^{ \imath \theta Z})\, , \quad \mbox{for every $\theta \in \R$}.\]

For sums of power of divisors of integers, we use the following notations.  For any integer $c \ge 0$, we denote with $\sigma_c(m)$ the sum of the $c$th powers of the divisors of $m$:
\[
\sigma_c(m)=\sum_{j\mid m} j^c\,, \quad \mbox{for $m \ge 1$}\,.
\]
Moreover, $\sigma_c^{\text{odd}}(m)$ denotes the analogous sum but restricted to the odd divisors of~$m$:
\[\sigma_c^{\text{odd}}(m)=\sum_{j\mid m,\, j\,\mbox{\tiny odd}} j^c\,,\quad \mbox{for $m \ge 1$}\,.\]

\subsection{Plan of the paper} Section \ref{seccion:Khinchin families} covers the basic background material on Khinchin families, the framework of the present paper, while Section \ref{seccion:gaussian Khinchin families} focuses on Gaussian and strongly Gaussian Khinchin families and the Hayman class. The reference \cite{K_uno} covers Khinchin  families at length.

Section \ref{seccion:set construction} describes the fundamental set constructions from the point of view of Khinchin families and the exponential function, which furnishes the basic context of application of the results of this paper.

The basic criteria for verifying that the exponential of a power series of nonnegative 
 coefficients is in the Hayman class are the main results of Section \ref{seccion:exponentials and Hayman}.

Finally, in Section \ref{seccion:asymptotic formulas coefficients} we briefly recall the procedure to obtain asymptotic  formulas of coefficients of strongly Gaussian power series.

\section{Khinchin families}
\label{seccion:Khinchin families} We denote by $\K$ the class of \textit{nonconstant} power series $f(z)=\sum_{n=0}^\infty a_n z^n$  with positive radius of convergence,  which have  \textit{nonnegative} Taylor coefficients, $a_n \ge 0$, for each $n \ge 0$, and such that  $a_0>0$.

The \textit{Khinchin  family} of  such a  power series $f \in \K$ with radius of convergence $R>0$ is the family of random variables $(X_t)_{t \in [0,R)}$ with values in $\{0, 1, \ldots\}$ and with mass functions given by
\[
\P(X_t=n)=\frac{a_n t^n}{f(t)}\, , \quad \mbox{for each $n \ge 0$ and $t\in (0,R)$}.
\]
The variable $X_0$ of the family is defined as $X_0\equiv 0$. Notice that $f(t)>0$ for each $t \in[0,R)$.

Any Khinchin family is continuous in distribution in $[0,R)$, in the sense that if a sequence $(s_n)_{n \ge 1} \subset [0,R)$ converges to $s_0\in [0,R)$, then $X_{s_n}\stackrel{\rm d}{\longrightarrow}
X_{s_0}$, as $n \to \infty$ (see, for instance, \cite{Kq_uno}). Observe that no hypothesis upon joint distribution of the variables~$X_t$ is considered; $(X_t)_{t \in [0,R)}$ is a family, not a process.

\subsection{Basic properties}
For the basic theory of Khinchin families (results, proofs, examples and applications), we refer the reader to \cite{K_uno}. Here we describe  the specific aspects of the theory to be used in the present  paper.

\subsubsection{Mean and variance functions}
For the mean and variance of $X_t$ we reserve the notation
$m(t)=\E(X_t)$ and $\sigma^2(t)=\V(X_t)$, for $t \in [0,R)$. In terms of $f$, the  mean and the variance of $X_t$ may be  written as
\[m(t)=\frac{t f^\prime(t)}{f(t)}, \qquad \sigma^2(t)=t m^\prime(t)\, , \quad \mbox{for $t \in [0,R)$}\,.\]

For each $t \in (0,R)$, the variable $X_t$ is not a constant, and so $\sigma^2(t)>0$. Consequently,  $m(t)$ is strictly increasing in $[0,R)$, though, in general, $\sigma(t)$ is not increasing. We  denote \[M_f=\lim_{t \uparrow R} m(t)\,.\]

Whenever $M_f=+\infty$, for each integer $n \ge 0$, we use $t_n$ to denote the unique $t_n\in [0,R)$, such that $m(t_n)=n$.

\subsubsection{Normalization and characteristic functions}

For each $t \in (0,R)$, the normalization of $X_t$ is  
\begin{equation*}
\breve{X}_t\triangleq\frac{X_t-m(t)}{\sigma(t)}\,\cdot
\end{equation*}
The characteristic function of the variable $X_t$ may be written in terms of the power series $f$ itself as
\begin{equation*}
\E(e^{\imath \theta X_t})=\frac{f(te^{\imath \theta})}{f(t)}\, , \quad \mbox{for $t\in (0,R)$ and $\theta \in \R$}\,,
\end{equation*}
while for its normalized version $\breve{X}_t$ we have that
\begin{equation*}
\E(e^{\imath \theta \breve{X}_t})=\E(e^{\imath \theta X_t/\sigma(t)}) \,e^{-\imath \theta m(t)/\sigma(t)}\, , \quad \mbox{for $t\in (0,R)$ and $\theta \in \R$}\,,
\end{equation*}
and so,
\[
|\E(e^{\imath \theta \breve{X}_t})|=|\E(e^{\imath \theta X_t/\sigma(t)})|\, , \quad \mbox{for $t\in (0,R)$ and  $\theta \in \R$}.
\]

\subsubsection{Fulcrum $F$ of $f$}\label{section:auxiliary F} The holomorphic function $f$ does not vanish on the real interval $[0,R)$, and so, it does not vanish in a simply connected region containing that interval.  We may consider $\ln f$, a holomorphic branch of the logarithm of $f$ which is real on $[0,R)$, and the  function $F$, which we shall call
 the \textit{fulcrum of $f$},  defined and holomorphic in a region containing $(-\infty, \ln R)$, and which is given by
\[F(z)=\ln f(e^z).\]

If $f$ does not vanish anywhere in the disk $\D(0,R)$, then the fulcrum $F$ of $f$ is defined in the whole half plane $\{z\in \C: \Re z< \ln R\}$.
In this paper,  this situation of~$f$ nonvanishing in $\D(0,R)$ is the most interesting. In this case $f(z)=e^{g(z)}$, where~$g$ is a function holomorphic in $\D(0,R)$ and $g(0)\in \R$, and the fulcrum $F$ of $f$ may be written~as
\[
F(z)=g(e^z)\,, \quad \mbox{for $z$ such that $\Re z<\ln R$}\,.
\]
The mean and variance function of $f$ may be expressed  in terms of its fulcrum $F$ as
\[
m(t)=F^\prime(s)\quad \mbox{and} \quad \sigma^2(t)=F^{\prime\prime}(s)\,,  \quad \mbox{for $s <\ln R$ and $t=e^s$}\,.
\]

\subsection{Hayman's identity}

For a power series $f(z)=\sum_{n=0}^\infty a_n z^n$ in $\K$, Cauchy's formula for the coefficient $a_n$ may be written in terms of the characteristic function of its Khinchin  family $(X_t)_{t \in [0,R)}$ as
\[a_n=\frac{f(t)}{2\pi t^n }\int_{|\theta|<\pi} \E(e^{\imath \theta X_t})\, e^{-\imath \theta n} \, d \theta\, , \quad \mbox{for each $t \in (0,R)$ and $n\ge 1$}.\]

In terms of the characteristic function of the normalized variable $\breve{X}_t$, Cauchy's formula becomes
\begin{equation*}
a_n=\frac{f(t)}{2\pi \,t^n \,\sigma(t)}\int\nolimits_{|\theta|<\pi\sigma(t)} \E(e^{\imath \theta \breve{X}_t}) \, e^{-\imath \theta (n-m(t))/\sigma(t)} \, d \theta\quad \mbox{for  $t \in (0,R)$ and $n \ge 1$}.
\end{equation*}

If $M_f=\infty$, we may take for each $n\ge 1$ the (unique) radius $t_n \in (0,R)$  so that $m(t_n)=n$, to write
\begin{equation}
\label{eq:haymans formula}
a_n=\frac{f(t_n)}{2\pi \,t_n^n \,\sigma(t_n)}\int_{|\theta|<\pi\sigma(t_n)} \E(e^{\imath \theta \breve{X}_{t_n}}) \, d \theta\, , \quad \mbox{for each $n \ge 1$}\, ,
\end{equation}
which we call \textit{Hayman's identity}.

This identity \eqref{eq:haymans formula}, which is just Cauchy's formula with an appropriate  choice of radius $t_n$, neatly encapsulates, in fact,  the saddle point method.

\subsection{Basic Khinchin families}\label{subsec:basic families} The most basic collections of probability distributions in $\{0,1, \ldots\}$, i.e., Bernoulli and binomial, geometric and negative binomial and Poisson, are (the most basic) Khinchin families.
A quick review follows; see more details, for example, in Section 2.1.6 of~\cite{K_uno}.
\begin{enumerate}[(a)]
\item The  Khinchin family of the function $f(z)=1+z$ consists of the \textit{Bernoulli variables}.

In this case, $R=\infty$, and the mean and variance functions are $m(t)=t/(1+t)$ and $\sigma^2(t)=t/(1+t)^2$. For each $t >0$, the random variable
 $X_t$ is a \textit{Bernoulli}  variable  with parameter $p=t/(1+t)$.

\item Let $f(z)=(1+z)^N$, with integer $N \ge 1$. This is the \textit{binomial case}.

In this case, $R=\infty$, and the mean and variance functions are $m(t)=Nt/(1+t)$ and $\sigma^2(t)=Nt/(1+t)^2$. For each $t >0$, the random variable
 $X_t$ is a
  \textit{binomial}  variable
   with parameters $N$ and $p=t/(1+t)$.

 \item The function $f(z)=1/(1-z)$ corresponds to the \textit{geometric case}.

 In this case  $R=1$, and the mean and variance functions are $m(t)=t/(1-t)$ and $\sigma^2(t)=t/(1-t)^2$. For each $t\in (0,1)$, the random variable $X_t$ is a \textit{geometric} 
  variable (number of failures until first success) of parameter $1-t$, that is,
$\P(X_t=k)=t^{k} (1-t)$ for $k \ge 0$. 

\item  Let $f(z)=1/(1-z)^N$, with integer $N \ge 1$; the \textit{negative binomial} case.

In this case  $R=1$, and the mean and variance functions are $m(t)=N t/(1-t)$ and $\sigma^2(t)=N t/(1-t)^2$.
For each $t\in (0,1)$, the random variable $X_t$ is a \textit{negative binomial} variable  of parameters $N \ge 1$ and  $p=1-t$.

\item  The Khinchin family of the exponential function $f(z)=e^z$ consists of the Poisson variables: the \textit{Poisson case}.

In this case  $R=\infty$, and the mean and variance functions are $m(t)=t$ and $\sigma^2(t)=t$. For each $t>0$, the random variable $X_t$ in its Khinchin family is a \textit{Poisson} variable  with parameter $t$.

\end{enumerate}

\subsection{A couple of power series comparisons}\label{seccion:moments of basic families}
In the proofs of Theorems  \ref{teor:general for quasiexponentials} and~\ref{teor:general for quasigeometric}, the main results of this paper, we will resort to the following asymptotics for a couple of series.

\begin{propo}\label{prop:binomicos} For $\beta>0$,
\[\sum_{n=1}^\infty n^{\beta-1}\, t^n \sim \Gamma(\beta)\,\frac{1}{(1-t)^{\beta}}\,, \quad \mbox{as $t\uparrow 1$}\,.\]
\end{propo}
\begin{proof}
The binomial expansion gives that
\[\frac{1}{(1-z)^\beta}=\sum_{n=0}^\infty \frac{\Gamma(n+\beta)}{\Gamma(\beta)\, n!}\, z^n\,,\quad \mbox{for $|z|<1$}\,.\]
And Stirling's formula gives, for $\beta>0$,  that
\begin{equation}\label{eq:comparison coeffs of (1-z)beta}
\frac{\Gamma(n+\beta)}{\Gamma(\beta)\, n!}\sim \frac{n^{\beta-1}}{\Gamma(\beta)}\,, \quad \mbox{as $n \to \infty$}\,.\end{equation}
Thus,
\[
\sum_{n=1}^\infty n^{\beta-1} \, t^n \sim \Gamma(\beta)\, \frac{1}{(1-t)^{\beta}}\,, \quad \mbox{as $t\uparrow 1$}\,.\qedhere\]
\end{proof}

\begin{remark}[Moments of geometric random variables] The power series comparison of Proposition \ref{prop:binomicos} translates into a comparison of moments of the Khinchin family $(X_t)_{t \in [0,1)}$ of $1/(1-z)$, i.e., of geometric variables. Namely, the following:
for $\beta>0$,
\[
\E(X_t^\beta)\sim \Gamma(\beta+1)\,\E(X_t)^\beta\,, \quad \mbox{as $t \uparrow 1$}\,.
\]
For $\beta \in (-1,0)$, the same result holds for the variables  $
\tilde X_t$, the $X_t$ \emph{conditioned at being positive}, given by
$
\P(\tilde X_t=k)=t^{k{-}1} (1-t)$ for $k \ge 1$ and $t \in (0,1)$.
\end{remark}

\begin{propo}\label{prop:binomicos factoriales} For $\beta\ge 0$, we have that
\[
\sum_{n=0}^\infty n^\beta \,\frac{t^n}{n!}\sim t^\beta e^t 
\quad\text{and}\quad
\sum_{n=1}^\infty \frac{1}{n^\beta} \,\frac{t^n}{n!}\sim \frac{e^t}{t^\beta} \quad \mbox{as $t \to \infty$}.
\]
\end{propo}

For a proof of Proposition \ref{prop:binomicos factoriales}, we refer the reader to \cite{K_tres}.

\begin{remark}[Moments of Poisson random variables]
Proposition \ref{prop:binomicos factoriales} is actually a statement about positive moments of the Poisson variables $X_t$, and also about negative moments of the conditioned  Poisson variables $\tilde  X_t$, as the mean $t$ tends to $\infty$. The conditioned Poisson variables $\tilde X_t$ are given by
\[\P(\tilde X_t=k)=\frac{e^{-t}}{1-e^{-t}} \,\frac{t^k}{k!}\, , \quad \mbox{for $k \ge 1$ and $t >0$}.\]
This moments estimation is the following: if\/ $(X_t)_{t \in [0,1)}$ is the Khinchin family of\/ $e^z$, then
for any $\beta \ge 0$,
\[\E(X_t^\beta)\sim t^\beta\,, \quad \mbox{as $t \to \infty$},\]
and for  any $\beta<0$,
\[\E\big(\tilde X_t^\beta\big)\sim t^\beta\,, \quad \mbox{as $t \to \infty$}.\]
 See also \cite{JonesZ} for precise asymptotic \textit{expansions} of the negative moments of Poisson variables.
\end{remark}

\begin{remark}[Khinchin families and clans]
We refer the reader to \cite{K_tres} for more information on the asymptotic comparison of $\E(X_t^\beta)$ with $\E(X_t)^\beta$, as $t \uparrow R$, for a given Khinchin family $(X_t)$, which is behind the notion of a \emph{Khinchin clan}.
\end{remark}

\section{Set construction and Khinchin families}\label{seccion:set construction}

The \textit{set construction of combinatorial classes} fits nicely within  the framework of Khinchin  families. As a general reference for combinatorial classes and operations with them, we strongly suggest Chapter II in Flajolet--Sedgewick~\cite{Flajolet}. See also \cite{ArriataBarbourTavare}, where the constructions below are \textit{presented}  as \textit{decomposable structures}.

We shall abbreviate `exponential generating function' by egf, and `ordinary  generating function'   by ogf.

Most of the (ordinary or exponential) generating functions 
 of sets of combinatorial classes are exponentials of power series with nonnegative coefficients, \emph{the object of interest} of this paper.

\smallskip

We will verify by means of Theorems \ref{teor:general for quasiexponentials} and \ref{teor:general for quasigeometric} that most of the set constructions give rise to generating functions which are in the Hayman class, see Section \ref{section:Hayman class}, so that there are asymptotic formulas for their coefficients given by the Hayman and B\'{a}ez-Duarte asymptotic formulas of
Section \ref{seccion:formula for coefficients}.

\subsection{Labeled combinatorial classes and sets}
If $g(z)=\sum_{n=1}^\infty (b_n/n!) \,z^n$ is the exponential generating function (egf) of a labeled combinatorial class $\mathcal{G}$ (no object of size 0), then $f(z)=e^{g(z)}=\sum_{n=0}^{\infty}(a_n/n!) z^n$  is the egf of the labeled class of sets formed with the objects of $\mathcal{G}$, termed \textit{assemblies} in \cite{ArriataBarbourTavare}.

\smallskip

Next, we consider first sets of the basic classes: sets, lists and cycles, and then sets of rooted trees and of functions.

\subsubsection{Sets of sets} 

Let $\mathcal{G}$ be the labeled class of nonempty sets: $b_n=1$, for each $n \ge 1$, and $b_0=0$. Its egf is $g(z)=e^z-1$. Then
\[ f(z)=\exp(e^z-1)=\sum_{n=0}^\infty \frac{B_n}{n!}\, z^n
\]
is the egf of sets of sets, or, equivalently, of partitions of sets, see~\cite{Flajolet}, p. 107. Here, $B_n$ is the $n$th Bell number, which counts the number of partitions of the set $\{1, \ldots, n\}$.  In this case $R=\infty$,  and the mean and variance function are $m(t)=t e^t$ and $\sigma^2(t)=t(t+1)e^t$.


The characteristic function of the Khinchin family of $f$ is given by
\[\E(e^{\imath \theta X_t})=\exp (e^{t e^{\imath \theta}}-e^t )\, , \quad \mbox{for  $\theta \in \R$ and $t >0$}\, ,\]
and thus, for  $\theta \in \R$ and $t >0$,
\begin{equation}\label{formula de func caracteristica de bell}\E(e^{\imath\theta \breve{X}_t})=
\exp\big(e^{t e^{\imath\theta  {e^{-t/2}}/{\sqrt{t(t+1)}}}}-e^t-\imath \theta \sqrt{t/(t+1)}e^{t/2}\big)
.\end{equation}

\smallskip

$\bullet$ The class $\mathcal{G}$ of pointed sets (i.e., sets with one marked element) has egf $g(z)=ze^z$. The class of sets of pointed sets is isomorphic to the class of idempotent maps; its egf~$f$ is given by $f(z)=\exp(ze^z)$. See \cite{Flajolet}, p. 131.

\subsubsection{Sets of lists} For the labeled class $\mathcal{G}$ of (nonempty) lists, the  function $g$ is just $g(z)=z/(1-z)$. And the function $f(z)=\exp(z/(1-z))$ is the egf of the sets of lists, the so called fragmented permutations. See \cite{Flajolet}, p. 125.

\subsubsection{Sets of cycles} The labeled class $\mathcal{G}$ of (nonempty) cycles has $g(z)=\ln (1/(1-z))$ as egf. The function $f(z)=\exp(g(z))={1}/(1-z)$ is the egf of the sets of cycles,  or, in other terms, of  the permutations.

\smallskip
  $\bullet$ The length of the cycles could be restricted. Thus, for integer $k \ge 1$,
\[
\exp\Big(\sum_{n\le k} z^n/n\Big)\quad\text{or}\quad  \exp\Big(\sum_{n\ge k} z^n/n\Big)
\]
are  the egfs of permutations such that all the cycles in their cycle decomposition have length at most $k$, or 
 at least $k$, respectively. We may also consider, for~$k \ge 1$,
\[
\exp\Big(\sum_{d\ge 1,\, d\mid k} z^d/d\Big),
\]
which is the egf of the permutations $\sigma$ such that $\sigma^k$ is the identity.

\subsubsection{Sets of trees and sets of functions}\mbox{}\vspace{6pt}

  $\bullet$ \emph{Sets of trees \textup{(}forests\textup{)}}. The class $\mathcal{G}$ of  rooted (labeled) trees has egf $g(z)=\sum_{n=1}^\infty (n^{n-1}/n!) \,z^n$. See \cite{Flajolet}, Section II.5.1; this is Cayley's theorem. The class of forests (sets) of rooted (labeled) trees has then egf $f=e^g$.

Cayley's theorem also shows that the egf of the class of  unrooted (labeled) trees is $g(z)=\sum_{n=1}^\infty (n^{n-2}/n!)\, z^n$.

\smallskip

  $\bullet$ \emph{Sets of functions}. The class $\mathcal{G}$ of functions has egf $g(z)=\sum_{n=1}^\infty (n^n/n!) z^n$. See~\cite{Flajolet}, Section II.5.2. The class of sets of functions has then egf $f=e^g$.

\subsection{Unlabeled combinatorial classes and sets}

We split the discussion into multisets and (proper) sets. See Section I.2.2 in \cite{Flajolet} and  \cite{ArriataBarbourTavare} as general references for the set constructions of unlabeled combinatorial classes.

\subsubsection{Multisets of unlabeled combinatorial classes} If $C(z)=\sum_{n=1}^\infty c_n z^n$ is the ordinary generating function (ogf) of a combinatorial (unlabeled) class $\mathcal{G}$ (no object of size 0), then
\[
f(z)=\prod_{j=1}^\infty \frac{1}{(1-z^j)^{c_j}}=\exp\Big(\sum_{n=1}^\infty C(z^n)/n\Big)
\]
is the ogf  of the  \textit{class of sets formed with the objects of\/ $\mathcal{G}$}, termed  \textit{multisets} in \cite{ArriataBarbourTavare}.

\smallskip

We may write $f(z)=e^{g(z)}$, where $g$ is the power series
\[g(z)=\sum_{j,k\ge 1} c_j\,\frac{z^{kj}}{k}=\sum_{m=1}^\infty\Big(\sum_{j \mid m}j c_j\Big)\frac{z^m}{m}\triangleq \sum_{m=1}^\infty b_m\, z^m\,.\]
Observe that the power series $g$ has nonnegative coefficients:
\[
b_m=\frac{1}{m}\sum_{j \mid m}j c_j\,, \quad \mbox{for $m\ge 1$}\,.\]

\medskip

$\bullet$  The ogf of partitions, \textit{the partition function}, given by
\[P(z)=\prod_{j=1}^\infty \frac{1}{1-z^j}=\sum_{n=0}^\infty p(n)\, z^n\, , \quad \mbox{for $z \in \D$}\,,\]
is in $\K$. Here $p(n)$ is the number of partitions of the integer $n \ge 1$. 

The mean and variance functions
 of its Khinchin family are not so direct.  In this instance, $C$ is $C(z)=z/(1-z)$, as $c_n=1$ for $n \ge 1$ (one object of weight $n$, for each $n \ge 1$). And thus $g$ has coefficients $b_m=\sigma_1(m)/m$, where $\sigma_1(m)$ is  the sum of the divisors of the integer $m$.

\smallskip

For general $C(z)=\sum_{n=1}^\infty c_n z^n$  as above, the corresponding $f$ is the ogf of the colored partitions; $c_j$ different colors for part $j$.

\subsubsection{Sets of unlabeled combinatorial classes} If again we write $C(z)=\sum_{n=1}^\infty c_n z^n$ for the ordinary generating function (ogf) of a combinatorial (unlabeled) class $\mathcal{G}$ (no object of size 0), then \[f(z)=\prod_{j=1}^\infty (1+z^j)^{c_j}=\exp\Big(\sum_{n=1}^\infty (-1)^{n+1}C(z^n)/n\Big),\]
is the ogf  of the  class of (proper) sets formed with the objects of $\mathcal{G}$, termed \textit{selections} in \cite{ArriataBarbourTavare}.

\smallskip

We may write $f(z)=e^{g(z)}$, where $g$ is the power series
\[
g(z)=\sum_{j,k\ge 1} c_j\,\frac{z^{kj}\,(-1)^{k+1}}{k}=\sum_{m=1}^\infty\Big(\sum_{jk=m}j c_j\,(-1)^{k+1}\Big)\frac{z^m}{m}\,\cdot
\]

In general, the power series $g$ could have negative coefficients; this is the case, for instance, for $f(z)=(1+z)^5 (1+z^2)$.

But for the sequence of coefficients of $C$ given by $c_j=j^{c-1}, j \ge1$,  where $c$ is an integer $c\ge 1$, \textit{the coefficients of $g$ are nonnegative}; in fact, the $m$th coefficient of $g$ is
\[
\frac{1}{m} \,\sigma_c^{\text{odd}}(m) \,\omega(m)\,,
\]
where $\omega(m)$ is given by
\[\omega(m)=\frac{2^c-2}{2^c-1}\, 2^{\chi(m)c}+\frac{1}{2^c-1}\,, \quad \mbox{for $m \ge 1$}\,,\]
$\chi(m)$ is the highest integer exponent  so that $2^{\chi(m)}\mid m$, and $\sigma_c^{\text{odd}}(m)$ is the sum of the~$c$th powers of the odd divisors of $m$.
This is so because of the identity
\begin{equation} \label{eq:alternating sums of powers of divisors}\sum_{jk=m}j^c (-1)^{k+1}=\Big(\sum_{j\mid m,\, j\, \mbox{\tiny odd}} j^c\Big) \, \omega(m)\,, \quad \mbox{for $m \ge 1$}\,.\end{equation}

To verify \eqref{eq:alternating sums of powers of divisors},  observe first that as functions of $m$, both summations in \eqref{eq:alternating sums of powers of divisors}  are multiplicative. Write $m$ as $m=2^{\chi(m)} s$, with $s$ odd, and observe that  for $m=2^r$, with $r \ge 1$, the summation on the left is $\omega(m)$ and the summation on the right is 1, while for $m=s$ odd, the two sums coincide and $\omega(s)=1$.


We may also write
\[\sum_{jk=m}j^c (-1)^{k+1}=\sigma_c(m)-2 \sigma_c(m/2)\,, \quad \mbox{for $m \ge 1$}\,,
\]
with the understanding that $\sigma_c(m/2)=0$ if $m$ is odd. For the identity \eqref{eq:alternating sums of powers of divisors}  and a variety of relations among a number of diverse sums on divisors, we refer to \cite{Glaisher}.

More generally,  if the coefficient $c_j$ of $C$ is given by $c_j=R(j)$, for $j \ge 1$, where $R(z)$ is a polynomial with nonnegative integer coefficients, then the coefficients of $g$ are nonnegative.

\smallskip

$\bullet$  The ogf of \textit{partitions with distinct parts}, given by
\[Q(z)=\prod_{j=1}^\infty (1+z^j)=\sum_{n=0}^\infty q(n)\, z^n\, , \quad \mbox{for $z \in \D$},\]
is in $\K$. Here $q(n)$ is the number of partitions into distinct parts of the integer $n \ge 1$.

In this particular instance, the power series $C$ is $C(z)=z/(1-z)$, as  $c_n=1$, for $n \ge 1$, and the power series $g$ is
\[
g(z)=\sum_{m=1}^\infty \frac{\sigma_1(m)-2\sigma_1(m/2)}{m}\, z^m\,.\]

\section{Gaussian Khinchin families and Hayman class}\label{seccion:gaussian Khinchin families}

\begin{defin}
A power series $f \in \K$ and its  Khinchin family $(X_t)_{t \in [0,R)}$ are termed \textit{Gaussian} if
the normalized sequence $(\breve{X}_t)$ converges in distribution, as $t\uparrow R$, to the standard normal or, equivalently, if
\[\lim_{t \uparrow R} \E(e^{\imath \theta \breve{X}_t})=e^{-\theta^2/2}\, , \quad \mbox{for each $\theta \in \R$}.\]
\end{defin}

Among the
basic Khinchin families considered in Section~\ref{subsec:basic families},
only the family associated to the exponential $e^z$ is Gaussian.

For the functions $f(z)=1+z$ (Bernoulli case) or $f(z)=(1+z)^N$ (binomial case), the corresponding $(\breve{X}_t)$ converges in distribution, as $t\to\infty $, to the constant~0. In fact, for any polynomial in $\K $, the corresponding $(\breve{X}_t)$ converges in distribution, as $t\to\infty$, towards the constant~0.

For the function $f(z)=1/(1-z)$ (geometric case),
 $(\breve{X}_t)$ converges in distribution, as $t\uparrow 1$, towards a variable $Z$, where  $Z+1$ is an exponential variable  of parameter~1. For $f(z)=1/(1-z)^N$, with integer $N\ge 1$ (negative binomial case),
it can be verified that
 $(\breve{X}_t)$ converges in distribution, as $t\uparrow 1$, towards a variable $Z_N$, where  $Z_N+\sqrt{N}$ follows a Gamma distribution with shape parameter $N$ and rate parameter $\sqrt{N}$ (or scale parameter $1/\sqrt{N}$).

Thus, in all the cases considered above, $(\breve{X}_t)$ converges in distribution as $t \uparrow R$, but only for the exponential $e^z$ the limit is the (standard) normal distribution.  See details in Section 3.1 of \cite{K_uno}.

\subsection{Gaussianity of exponentials}\label{section:gaussianity exponentials} The following Theorem \ref{teor:hayman} is  the  basic criterion for gaussianity in terms of the fulcrum $F$ of $f$
(of Section \ref{section:auxiliary F}); it originates in Hayman's \cite{Hayman}, but we refer to  Theorem 3.2 in \cite{K_uno} for a proof.

\begin{theorem}\label{teor:hayman}  If $f\in \K$ has radius of convergence $R>0$ and vanishes nowhere in $\D(0,R)$, and if for the fulcrum  $F$ of $f$ one has
\begin{equation}\label{eq:condicion de Hayman}\lim_{s \uparrow \ln R} \frac{\sup_{\phi\in\R} \big|F^{\prime\prime\prime}(s+i\phi)\big|}{F^{\prime\prime}(s)^{3/2}}=0\, ,\end{equation}
then $f$ is Gaussian.
\end{theorem}

Notice that Theorem \ref{teor:hayman} applies whenever $f \in \K$ is of the form $f=e^g$ for some power series $g$, that may have negative Taylor coefficients. As registered in the following Theorem \ref{teor:hayman bis}, if $f=e^g$ and $g$ has nonnegative coefficients, a simpler condition on $g$ implies the gaussianity of $f=e^g$.

\begin{theorem}\label{teor:hayman bis}
Let $f\in \K $ be such that $f=e^g$, where $g$ has radius of convergence $R>0$ and nonnegative coefficients. If~$g$ is a polynomial of degree $1$ or $g$ satisfies
\begin{equation}
\label{eq:condition on g for f=exp(g) Gaussian}
\lim_{t \uparrow R} \frac{g^{\prime\prime\prime}(t)}{{g^{\prime\prime}(t)}^{3/2}}=0\, ,
\end{equation}
then $f$ is Gaussian. \end{theorem}

For a proof, see Theorem 3.3 in \cite{K_uno}.

\medskip

$\bullet$ For the exponential  function $f(z)=e^z$, we have that its fulcrum is $F(z)=e^z$, and its gaussianity follows readily from Theorem \ref{teor:hayman}.

\vskip2pt

$\bullet$ More generally, if $B(z)=\sum_{j=0}^N b_j \, z^j$  is a polynomial of degree $N$ such that  $e^{B}\in \K $, then  $e^{B}$ is Gaussian. For, in this case, $F(z)=B(e^z)$, and for $z=s+\imath \phi$ we have that
\[|F^{\prime\prime \prime}(z)|=\Big|\sum_{j=0}^N b_j \,j^3 e^{jz}\Big| \le \sum_{j=0}^N |b_j| \,j^3 e^{js}=O(e^{Ns})\, ,\]
 while
 \begin{equation}
 \label{eq:second derivative fulcrum poly}
 F^{\prime\prime}(s)=\sum_{j=0}^N b_j \,j^2 e^{js}\sim b_N \,N^2 e^{Ns}\, \quad \mbox{as $s \uparrow \infty$}\, ,
 \end{equation}
 and gaussianity follows from Theorem \ref{teor:hayman}. Observe that \eqref{eq:second derivative fulcrum poly}  implies that $b_N$ is real and $b_N >0$; besides, since $e^{B} \in \K $, the coefficients of the polynomial $B$ must be real numbers.

\medskip

$\bullet$ For $f(z)=\exp(e^z-1)$, the egf of the \textit{class of sets of sets},  and for $f(z)=\exp(z/(1-z))$, the egf of the \textit{class of sets of  lists}, Theorem \ref{teor:hayman} (or Theorem \ref{teor:hayman bis})  gives readily that both are  Gaussian.
Similarly, the function $f(z)=\exp(z e^z)$, the egf of the \textit{class of sets of pointed sets}, is Gaussian.

\medskip

 $\bullet $ Again directly from  Theorem \ref{teor:hayman} or Theorem \ref{teor:hayman bis}, we see that the functions $f(z)=\exp(1/(1-z)^\gamma)$, with $\gamma>0$, are all Gaussian.   The functions
\[
f(z)=\exp\Big(\sum_{n=1}^\infty n^\alpha z^n\Big)\,, \quad \mbox{for $|z|<1$}\,,
\]
with $\alpha>-1$ are also Gaussian. Indeed, the fulcrum $F$ of such an $f$ is
\[F(z)=\sum_{n=1}^\infty n ^\alpha e^{n z}\, , \quad \mbox{for $\Re z<0$}\,.\]
Then, for $s<0$ and $r=e^s$, we have that
\[
\sup_{\phi \in \R}|F^{\prime\prime\prime}(s+\imath \phi)|=\sum_{n=1}^\infty n^{\alpha+3} \,r^n\,,\quad \mbox{while} \quad
|F^{\prime\prime}(s)|=\sum_{n=1}^\infty n^{\alpha+2} \,r^n\,.\]

By appealing to Proposition \ref{prop:binomicos}, 
 we deduce that condition \eqref{eq:condicion de Hayman} is satisfied. Observe that $R=1$, in this case.

\smallskip

  $\bullet$ In the same vein, the egf of \textit{sets of functions}
$f(z)=\exp\big(\sum_{n=1}^\infty ({n^n}/{n!}) z^n\big)$,
with $R=1/e$, is also Gaussian.
Similarly, and more generally,
\[
f(z)=\exp\Big(\sum_{n=1}^\infty \frac{n^{n-\delta}}{n!} \,z^n\Big)
\]
is Gaussian, for $0\le \delta <1/2$. To see this, observe (Stirling's formula) that
\[
\frac{n^{n-\delta}}{n!}\asymp \frac{e^n}{n^{\delta+1/2}}\,, \quad \mbox{as $n \to \infty$}\,.
\]
and argue as above with $\alpha=-\delta-1/2$.

\medskip

 $\bullet$ The \textit{partition function $P$ and the ogf $\,Q$ of partitions into distinct parts} are seen to be Gaussian as a consequence of Theorem \ref{teor:hayman}; but see Section \ref{seccion:exponentials and Hayman} for the stronger statement that $P$ and $Q$ are both in the Hayman class.

\medskip

Now we turn our attention to a couple of \textit{non-examples}.

\medskip

$\bullet$ The egf $f(z)=\exp(\ln (1/(1-z)))$ of the \textit{class of sets of  cycles}, which is simply $f(z)=1/(1-z)$, the egf of  permutations, is \textit{not Gaussian}, as we have mentioned at the beginning of this section.
  The condition of Theorem \ref{teor:hayman} is, of course, not satisfied; in fact, for the corresponding fulcrum $F(z)=-\ln (1-e^z)$, for $\Re z <0$, it holds that
\[
\lim_{s \uparrow 0} \frac{\sup_{\phi\in\R} |F^{\prime\prime\prime}(s+i\phi)|}{F^{\prime\prime}(s)^{3/2}}=2\,.
\]

Analogously,  the class of sets of cycles of length at least $k$ is not Gaussian, for any $k \ge 1$.
\medskip

 $\bullet$ The egf $f(z)=\exp\big(\sum_{n=1}^\infty ({n^{n-1}}/{n!}) z^n\big)$ of forests of rooted trees is $\textit{not Gaussian}$. In fact, we are going to check that the characteristic functions of its  normalized Khinchin family converge in distribution to the constant 1.

The power series $f$ has radius of convergence $R=1/e$, and
\[
\lim_{t\uparrow 1/e} f(t)=\exp\Big(\sum_{n=1}^\infty \frac{n^{n-1}}{n!\, e^n}\Big)<+\infty,
\]
since
\[
\frac{n^{n-1}}{n!\, e^n}\asymp \frac{1}{n^{3/2}}\quad\text{as $n\to\infty$.}
\]
In particular, the function $f$ extends to be continuous in $\mbox{\rm cl}(\D(0,1/e))$.

Now
\[
m(t)=\sum_{n=1}^\infty \frac{n^{n}}{n!} \,t^n, \quad \mbox{for $t \in (0,1/e)$}\,,
\]
and Proposition  \ref{prop:binomicos} gives that 
\[
m(t/e)=\sum_{n=1}^\infty \frac{n^{n}}{n!\,e^n} \,t^n \asymp \frac{1}{(1-t)^{1/2}}\,, \quad \mbox{as $t \uparrow 1$}\,.
\]
Also
\[
\sigma^2(t)=\sum_{n=1}^\infty \frac{n^{n+1}}{n!} \,t^n, \quad \mbox{for $t \in (0,1/e)$}\,,
\]
and Proposition  \ref{prop:binomicos} gives that
\[
\sigma^2(t/e)=\sum_{n=1}^\infty \frac{n^{n+1}}{n!\,e^n} \,t^n\asymp \frac{1}{(1-t)^{3/2}}\,, \quad \mbox{as $t \uparrow 1$}\,.
\]
Thus
\[
\frac{\sigma(t/e)}{m(t/e)}\asymp \frac{1}{(1-t)^{1/4}}\,, \quad \mbox{as $t \uparrow 1$}\,,\]
and so, in particular,
\begin{equation}\label{eq:m/q forests of rooted trees}
\lim_{t \uparrow 1/e} \frac{m(t)}{\sigma(t)}=0\,.\end{equation}

\medskip

Now for $\theta \in \R$ we may write
  \[
  \E(e^{\imath \theta \breve{X}_t})=\frac{f(t e^{\imath \theta /\sigma(t)})}{f(t)}\, e^{-\imath \theta m(t)/\sigma(t)}\,.
  \]
  The first factor of this expression of the characteristic function of $\breve{X}_t$ tends towards~1 because $\lim_{t \uparrow 1} \sigma(t)=+\infty$ and the continuity of $f$ on $\mbox{cl}(\D(0,1/e)$, while the second factor  tends towards 1,
  as a consequence of \eqref{eq:m/q forests of rooted trees}. And thus, $\lim_{t\uparrow 1/e}\E(e^{\imath \theta \breve{X}_t})=1$, for any $\theta \in \R$.

  This means that $\breve{X}_t$ tends in distribution towards the constant 0, and not to the Gaussian distribution.

\smallskip

Similarly, and more generally, for
\[
f(z)=\exp\Big(\sum_{n=1}^\infty \frac{n^{n-\delta}}{n!}\, z^n\Big)
\]
the corresponding $(\breve{X}_t)$ tends in distribution, as $t\uparrow 1$, towards the constant 0, for $1/2<\delta \le 5/2$, and thus,
\emph{the efg of the sets of functions is Gaussian, but the egfs of sets of trees and of sets of rooted trees are not Gaussian}.

\begin{question}
Is \[f(z)=\exp\Big(\sum_{n=1}^\infty \frac{n^{n-1/2}}{n!}\, z^n\Big)\] Gaussian? It is not strongly Gaussian (see the definition below)
 since for its Khinchin family we have that
\[\lim_{t\uparrow 1/e} \frac{m(t)}{\sigma(t)}=(2\pi)^{-1/4}\,,\] while for strongly Gaussian power series it is always the case that that limit is $+\infty$. See formula \eqref{eq:m/s a infty para strongly gaussian} below.
\end{question}

\subsection{Strongly Gaussian Khinchin families}

The notion of strongly Gaussian power series $f $ in $\K$ was introduced by B\'{a}ez-Duarte in  \cite{BaezDuarte}.

\begin{defin}\label{defin:strongly gaussian} A power series $f \in \K$ and its  Khinchin family $(X_t)_{t \in [0,R)}$ are termed  \textit{strongly Gaussian} if
\[\mathrm{a)} \quad\lim_{t \uparrow R} \sigma(t)=+\infty,\quad \mbox{and}\quad \mathrm{b)} \quad
\lim_{t \uparrow R} \int\nolimits_{|\theta|<\pi \sigma(t)}\big|\E(e^{\imath \theta \breve{X}_t})-e^{-\theta^2/2}\big| \, d\theta=0.\]\end{defin}

The exponential $f(z)=e^z$ is strongly Gaussian. In this case $\sigma(t)=\sqrt{t}$. This strong gaussianity follows from the gaussianity of $e^z$ and dominated convergence using the bound
\[
\big|\E(e^{\imath\theta \breve{X}_t})\big|=e^{t (\cos(\theta/\sqrt{t})-1)}\le e^{-2\theta^2/\pi^2}\, , \quad \mbox{for  $\theta\in \R,\, t>0$ such that  $|\theta|< \pi \sqrt{t}$}.\]

As we will see in Theorem \ref{teor:Hayman central limit}, strongly Gaussian power series are Gaussian. Thus the power series $(1+z)^N$ and $1/(1-z)^N$, for $N \ge 1$, are all not strongly Gaussian.

\smallskip

For strongly Gaussian power series, we have the following  key Theorem \ref{teor:Hayman central limit}. It appears in \cite{Hayman} for power series in the Hayman class. See Theorem~A in~\cite{K_uno} for a proof under  the more general asumption of strongly Gaussian power series.

\begin{theorem}[Hayman's   central limit theorem] \label{teor:Hayman central limit} If the power series $f(z)=\sum_{n=0}^\infty a_n z^n $ in $\K$ is strongly Gaussian, then
\begin{equation}\label{eq:Hayman central limit}\lim_{t \uparrow R} \sup\limits_{n \in \Z} \Big|\frac{a_n t^n}{f(t)} \sqrt{2\pi}\sigma(t)-e^{-(n-m(t))^2/(2\sigma^2(t))}\Big|=0\,.\end{equation}

Moreover,
\[
\lim_{t \uparrow R}\P( \breve{X}_t \le b)=\Phi(b)\,, \quad \mbox{for every $b \in \R$}\, ,\]
and so,   $(\breve{X}_t)$ converges in distribution towards the standard normal  and $f$ is Gaussian.
\end{theorem}
In this statement, $a_n=0$, for $n <0$.

\smallskip

By considering  $n={-}1$ in \eqref{eq:Hayman central limit}, it follows that
\begin{equation}
\label{eq:m/s a infty para strongly gaussian}
\lim_{t \uparrow R}\dfrac{m(t)}{\sigma(t)}=+\infty\,,
\end{equation}
and, in particular, that  $M_f=\infty$ for every  $f\in \K$ that is strongly Gaussian.

\smallskip

The function $f(z)=e^{z^2}$ is Gaussian, because of Theorem \ref{teor:hayman}. Its variance  function $\sigma^2(t)=4t^2$ tends towards $\infty$ as $t \to\infty$, but $f$ is  not strongly Gaussian, since its Taylor coefficients of odd order are null and  do not satisfy the asymptotic formula~\eqref{eq:asymptotic formula coeffs}  below.

\subsubsection{Coefficients of strongly Gaussian power series}\label{seccion:formula for coefficients}

For strongly Gaussian power series, we have the following asymptotic formula for its coefficients.

\begin{theorem}[Hayman's asymptotic formula]\label{teor:hayman asymptotic formula} If $f(z)=\sum_{n=0}^\infty a_n z^n$ in $\K$ is strongly Gaussian, then
\begin{equation}
\label{eq:asymptotic formula coeffs}
a_n \sim  \frac{1}{\sqrt{2\pi}} \,\frac{f(t_n)}{t_n^n \,\sigma(t_n)}\, , \quad \mbox{as $n \to \infty$} .
\end{equation}
\end{theorem}
In the asymptotic formula above,  $t_n$ is given by $m(t_n)=n$, for each $n \ge 1$. This asymptotic formula follows readily from Theorem \ref{teor:Hayman central limit}, or alternatively, from Hayman's formula \eqref{eq:haymans formula} and strong gaussianity.

For the exponential function $f(z)=e^z$, one has $m(t)=t$ and $\sigma(t)=\sqrt{t}$, for $t \ge 0$, and $t_n=n$ for $n \ge 1$. The asymptotic formula above gives
\[
\frac{1}{n!}\sim \frac{1}{\sqrt{2\pi}} \,\frac{e^n}{n^n \sqrt{n}}\, , \quad \mbox{as $n \to \infty$}\,,
\]
that is Stirling's formula.

Actually, if $\omega_n$ is a good approximation of $t_n$, in the sense that
\[\lim_{n \to \infty} \frac{m(\omega_n)-n}{\sigma(\omega_n)}=0\, ,\]
then from Hayman's formula \eqref{eq:haymans formula} and strong gaussianity we have that
\begin{equation}
\label{eq:formula asintotica coeficientes de Hayman bis}
a_n \sim \frac{1}{\sqrt{2\pi}}\,\frac{f(\omega_n)}{\sigma(\omega_n) \, \omega_n^n}\, ,\quad \mbox{as $n \to \infty$}\,.
\end{equation}

In general, precise expressions for the $t_n$ are rare, since inverting $m(t)$ is usually  complicated. But, fortunately, in practice, one can do with a certain  asymptotic approximation due to B\'{a}ez-Duarte, \cite{BaezDuarte},  which we now describe.

Suppose that $ f \in \K$ is strongly Gaussian. Assume that $\widetilde{m}(t)$ is continuous and monotonically increasing to $+\infty$ in $[0,R)$ and that $\widetilde{m}(t)$ is a good approximation of $m(t)$ in the sense that
\begin{equation}\label{eq:condition baez-duarte}
\lim_{t \uparrow R} \frac{m(t)-\widetilde{m}(t)}{\sigma(t)}=0\,.
\end{equation}
Let  $\tau_n$ be  defined by $\widetilde{m}(\tau_n)=n$, for each $n \ge 1$.

\begin{theorem}[B\'{a}ez-Duarte asymptotic formula]
\label{teor:baez-duarte}
With the notations above, if $f(z)=\sum_{n=0}^\infty a_n z^n$ in $\K$ is strongly Gaussian and \eqref{eq:condition baez-duarte} is satisfied, then
\[
a_n \sim  \frac{1}{\sqrt{2\pi}} \,\frac{f(\tau_n)}{\tau_n^n \,\sigma(\tau_n)}\, , \quad \mbox{as $n \to \infty$} .
\]
\end{theorem}

This follows readily from \eqref{eq:formula asintotica coeficientes de Hayman bis}.

Besides, if  $\widetilde{\sigma}(t)$ is such  that $\sigma(t) \sim  \widetilde{\sigma}(t)$ as $t \uparrow R$, we may further write
\begin{equation}
\label{eq:formula de hayman-baez-duarte}
a_n \sim  \frac{1}{\sqrt{2\pi}} \,\frac{f(\tau_n)}{\tau_n^n \,\widetilde{\sigma}(\tau_n)}\, , \quad \mbox{as $n \to \infty$} .
\end{equation}

In practice, using \eqref{eq:formula de hayman-baez-duarte} requires approximating $m$ by $\widetilde{m}$ to obtain $\tau_n$, and then obtaining good enough approximations of $\sigma$ and $f$ on $(0,R)$  to produce asymptotic formulas of $\sigma(\tau_n)$ and $f(\tau_n)$.

\subsection{Hayman class}\label{section:Hayman class}

The class of Hayman consists of power series $f$ in $\K$ which satisfy some concrete and verifiable conditions which imply that $f$ is strongly Gaussian, see Theorem \ref{teor:hayman implica baez} below.

\begin{defin}
A power series $f\in \K $ is in the \textit{Hayman class} {\upshape{(}}or is Hayman-admissible or just $H$-admissible{\upshape{)}} if for a certain function $h\colon  [0,R)\rightarrow (0, \pi]$, the following conditions are satisfied:
\begin{align}\label{eq:h mayor caracteristica}\mbox{\rm (major arc):} \qquad
&\lim_{t \uparrow R}\sup_{|\theta|\le h(t)\, \sigma(t)} \big|\E(e^{\imath \theta \breve{X}_t})\,e^{\theta^2/2}-1\big|=0\, , \\
\label{eq:h menor caracteristica}\mbox{ \rm (minor arc):} \qquad
&\lim_{t\uparrow R}\sigma(t) \, \sup_{h(t)\sigma(t)\le |\theta| \le \pi \sigma(t)} |\E(e^{\imath \theta \breve{X}_t})|=0\, ,
\intertext{and}
 \label{eq:condicion de varianza en Hayman}\mbox{\rm (variance condition):}\qquad &\lim_{t \uparrow R} \sigma(t)=\infty.\end{align}
 \end{defin}

We refer to the function $h$ in the definition above as a  \textit{cut} between a \textit{major} arc and a \textit{minor} arc.

Some authors include in the Hayman class power series with a finite number of negative coefficients. This is not the case in this paper.

 For $f$ in the Hayman class,  the characteristic  function of $\breve{X}_t$ is uniformly approximated by  $e^{-\theta^2/2}$ in  the major arc, while it is uniformly $o(1/\sigma(t))$ in the minor arc.

Observe that condition \eqref{eq:h menor caracteristica} may be written more simply and in terms of $f$ itself as the requirement that
\[
\lim_{t\uparrow R}\sigma(t) \, \sup_{h(t)\le |\theta| \le \pi } |\E(e^{\imath \theta X_t})|=\lim_{t\uparrow R}\sigma(t) \, \sup_{h(t)\le |\theta| \le \pi } \frac{|f(te^{\imath \theta})|}{f(t)}=0.
\]

\begin{theorem}\label{teor:hayman implica baez} Power series in the Hayman class are strongly Gaussian.\end{theorem}

In a certain sense, the above theorem places the conditions for Hayman class as a criterion for being strongly Gaussian. For a proof, see Theorem 3.7 in \cite{K_uno}. Thus we have that being in the Hayman class implies strong gaussianity, which in turn implies gaussianity. In applications, we shall always check that the power series belongs to the Hayman class.

We refer to \cite{K_uno},  and, of course,  to \cite{Hayman}, for examples and further properties of the Hayman class.

\section{Exponentials and the Hayman class}\label{seccion:exponentials and Hayman}

Along this section, we consider nonconstant power series $g(z)=\sum_{n=0}^\infty b_n z^n$ with \textit{nonnegative  coefficients} and radius of convergence $R>0$; it is not required that $g(0)>0$. We are interested in conditions on $g$ which guarantee that $f=e^g$, which is in $\K$, is in the Hayman class.

There are two results of Hayman in \cite{Hayman} along this line.
\begin{itemize}
\item[(a)] If $g$ is a power series in the Hayman class, then $f=e^g$ is in the Hayman class; see Theorem~VI in \cite{Hayman}.

    \item[(b)] If $B$ is a nonconstant polynomial with nonnegative coefficients and such that $Q_B=\gcd\{n \ge 1: b_n >0\}=1$, then $e^B$ is in the Hayman class.
\end{itemize}

\begin{remark}In \cite{OdlyzkoRichmond}, it is shown that if $g$ is in the Hayman class, then $f=e^g$ satisfies the stronger (than Hayman) conditions of Harris and Schoenfeld, \cite{HarrisSchoenfeld}, which allow to obtain full asymptotic expansions of coeficientes and not just asymptotic formulas.\end{remark}

\begin{remark}\label{remark:exponential poly hayman} In Theorem~X of \cite{Hayman}, Hayman proves the following stronger version of~(b): \textit{Let $B(z)=\sum_{n=0}^N b_n z^n$ be a nonconstant polynomial with real coefficients such that  for each $d>1$ there exists $m$, not a multiple of $d$, such that $b_m \neq 0$, and such that if $m(d)$ is the largest such $m$, then $b_{m(d)}>0$. If $e^B$ is in $\K$, then $e^B$ is in the Hayman class.}

For a direct proof of E2), we refer to Proposition 5.1 in \cite{K_uno}.
\end{remark}

Theorem \ref{teor:cut for e^g in Hayman class}, stated  below, which is Theorem 4.1 in \cite{K_uno},  exhibits conditions on the function $g$ which ensure that $f$ is in the Hayman class.

Further, we describe in Theorem \ref{teor:cut for e^g in Hayman class bis} a \lq practical\rq \, approach to verify that a power series satisfies the conditions of Theorem \ref{teor:cut for e^g in Hayman class}.

Finally, Theorems \ref{teor:general for quasiexponentials} and \ref{teor:general for quasigeometric} give  easily verifiable conditions on \textit{the coefficients~$b_n$ of~$g$} that imply that $f=e^g$ is in the Hayman class. These theorems  may be directly applied to the generating functions of sets constructions (and also to other exponentials) to check that they belong to the Hayman class.

\smallskip

Let $(X_t)_{t \in [0,R)}$ be the Khinchin family of  $f=e^g$. The mean and variance functions of~$f$, \textit{written in terms of~$g$}, are
\begin{align}
\label{eq:formula for in terms of g}
m(t)&=\frac{tf'(t)}{f(t)} = t g^\prime(t)=\sum_{n=1}^\infty n b_n \,t^n\,, \quad \mbox{for $t \in (0,R)$}\,,
\\
\intertext{and}   \label{eq:formula for sigma in terms of g}
\sigma^2(t)&=tm'(t)=tg^\prime(t)+t^2g^{\prime\prime}(t)=\sum_{n=1}^\infty n^2 b_n \,t^n\,,\quad \mbox{for $t \in (0,R)$}\,.
\end{align}
Since  $g$ has nonnegative  coefficients, the variance function  $\sigma^2(t)$ of $f$ is increasing in $[0,R)$. Observe also that $m(t)\le \sigma^2(t)$, for $t \in [0,R)$.

\smallskip

The variance condition \eqref{eq:condicion de varianza en Hayman} required for  $f$ to be in the Hayman class translates readily in the following condition in terms of $g$:
\begin{equation*}\lim_{t \uparrow R} \big(tg^\prime(t) +t^2g^{\prime\prime}(t)\big)=+\infty.\end{equation*}

\smallskip

To properly handle the minor and major arc conditions in terms of $g$, we introduce
\begin{equation}\label{definition omegag}
\omega_g(t)\triangleq\frac{1}{6}(b_1 t+8 b_2 t^2+\frac{9}{2} \,t^3g^{\prime\prime\prime}(t)) \, , \quad \mbox{for $t \in (0,R)$}.
\end{equation}

 \smallskip

The following theorem is Theorem 4.1 in \cite{K_uno}.
\begin{theorem}\label{teor:cut for e^g in Hayman class}
Let $g$ be a nonconstant power series with radius of convergence $R$ and nonnegative coefficients.

If the variance condition \begin{equation}\label{eq:variance condition for e^g}\lim_{t \uparrow R} \big(tg^\prime(t) +t^2g^{\prime\prime}(t)\big)=+\infty.\end{equation} is satisfied and there is a cut function $h\colon  [0,R) \to (0,\pi)$ so that \begin{equation}\label{eq:arco mayor para f=e^g} \lim_{t \uparrow R} \omega_g(t)\, h(t)^3=0\,,
\end{equation}
 and  \begin{equation}\label{eq:arco menor para f=e^g}
\lim_{t \uparrow R} \sigma(t) \exp\Big(\sup\limits_{h(t)\le |\theta|\le \pi}\Re g(te^{\imath \theta})-g(t)\Big)=0\end{equation}
 hold,   then $f=e^g$ is  in the Hayman class.
\end{theorem}

\medskip

Condition \eqref{eq:arco mayor para f=e^g} of Theorem \ref{teor:cut for e^g in Hayman class} gives that the cut function $h$  fulfills  condition \eqref{eq:h mayor caracteristica} on the major arc.
Also, condition  \eqref{eq:arco menor para f=e^g} of Theorem \ref{teor:cut for e^g in Hayman class}, which involves $h$ and $g$,  implies that condition \eqref{eq:h menor caracteristica} on the minor arc is satisfied.

\smallskip

Condition \eqref{eq:arco menor para f=e^g} on the minor arc is the most delicate  to check. In practice, it depends on properly bounding $\sup\nolimits_{\omega\le |\theta|\le \pi}\Re g(te^{\imath \theta})-g(t)$ for general $\omega$, as it is exhibited  in the following variant, actually a corollary, of Theorem \ref{teor:cut for e^g in Hayman class}.

\begin{theo}\label{teor:cut for e^g in Hayman class bis} Let $g$ be a nonconstant power series with radius of convergence $R$ and nonnegative coefficients, and let $h\colon  [0,R) \to (0,\pi)$ be a cut function.

Assume that the variance condition \eqref{eq:variance condition for e^g} and the condition on the major arc \eqref{eq:arco mayor para f=e^g} of Theorem {\rm\ref{teor:cut for e^g in Hayman class}} are satisfied.

Assume further that there are positive functions $U, V$ defined in $(t_0,R)$, for some  $t_0\in (0,R)$, where $U$ takes values in $(0,\pi]$ and $V$ in $(0,\infty)$, and  such that
\begin{equation}\label{eq:arco menor para f=e^g bis uno}
\sup\limits_{|\theta|\ge\omega}\big(\Re g(t e^{\imath \theta})-g(t)\big)\le -V(t) \, \omega^2\, , \quad \mbox{for $\omega\le U(t)$ and $t \in (t_0,R)$}\, ,
\end{equation}
and
\begin{equation}\label{eq:arco menor para f=e^g bis dos}
h(t)\le U(t), \quad \mbox{for $t \in (t_0,R)$}\quad \mbox{and} \quad  \lim_{t \uparrow R} \sigma(t) \, e^{-V(t)  h(t)^2}=0.
\end{equation}

Then $f=e^g$ is in the Hayman class.
\end{theo}
\begin{proof}
 The pair of conditions \eqref{eq:arco menor para f=e^g bis uno} and \eqref{eq:arco menor para f=e^g bis dos} together imply condition \eqref{eq:arco menor para f=e^g}. This is so since the first half of \eqref{eq:arco menor para f=e^g bis dos} allows us to take $\omega=h(t)$ in \eqref{eq:arco menor para f=e^g bis uno} to obtain
\[\exp\big(\sup\limits_{|\theta|\ge h(t)}\Re g(t e^{\imath \theta})-g(t)\big)\le \exp\big(-V(t) \, h(t)^2\big)\,, \quad \mbox{for $t\in (t_0,R)$}\,,\]
and to then apply the second half of \eqref{eq:arco menor para f=e^g bis dos}.\end{proof}

Notice further that if for a function $V$ defined in $(t_0,R)$, for some  $t_0\in (0,R)$, and taking values in $(0,+\infty)$,   we have that
\begin{equation}\label{eq:arco menor para f=e^g bis uno uno}
\Re g(t e^{\imath \theta})-g(t)\le -V(t) \, \theta^2\, , \quad \mbox{for $|\theta|\le \pi$ and $t \in (t_0,R)$}\, ,
\end{equation}
then we may set $U\equiv \pi$, and the cut $h$ is just required to satisfy
\begin{equation}\label{eq:arco menor para f=e^g bis dos dos}
\lim_{t \uparrow R} \sigma(t) \, e^{-V(t)  h(t)^2}=0.
\end{equation}
Thus these two conditions \eqref{eq:arco menor para f=e^g bis uno uno} and \eqref{eq:arco menor para f=e^g bis dos dos} imply the conditions \eqref{eq:arco menor para f=e^g bis uno} and \eqref{eq:arco menor para f=e^g bis dos} of Theorem \ref{teor:cut for e^g in Hayman class bis}.

\smallskip

Next we will exhibit easily verifiable  requirements  on the coefficients of $g$ to obtain functions $V$ and $U$ so that the conditions \eqref{eq:arco menor para f=e^g bis uno} and \eqref{eq:arco menor para f=e^g bis dos} hold, or, simply, a function~$V$ so that the conditions \eqref{eq:arco menor para f=e^g bis uno uno} and  \eqref{eq:arco menor para f=e^g bis dos dos} are satisfied.

The announced  requirements differ if the power series $g$ is entire or has finite radius of convergence; we split the discussion accordingly.

\subsection{Case \texorpdfstring{$g$}{g} entire, \texorpdfstring{$R=\infty$}{R=infinty}}

Here, $g$ is an entire  power series with nonnegative coefficients.

\begin{lem} \label{lema:major arc for quasiexponentials} If  an entire  power series $g$ with nonnegative coefficients satisfies, for some $\beta>0$ and $B>0$, that
\begin{equation}
\label{eq:major arc for quasiexponentials}
\Re g(te^{\imath \theta})-g(t)\le B(e^{\beta t \cos \theta }-e^{\beta t})\, , \quad \mbox{for $t >0$ and $|\theta|\le \pi$}\,,
\end{equation}
then condition \eqref{eq:arco menor para f=e^g bis uno uno} is satisfied with $V(t)=C e^{\beta t}$ for some constant $C>0$ depending only on $B$ and $\beta$,  and  condition \eqref{eq:arco menor para f=e^g bis dos dos} requires that
$\lim_{t \to \infty} \sigma(t) \exp(- C e^{\beta t} h(t)^2)=0$.
\end{lem}

\begin{proof}
For $t>0$ and $\beta>0$, the convexity of the exponential function $x \mapsto e^{\beta t x}$ in the interval $[-1,1]$ gives
\[
e^{\beta t}-e^{\beta t \cos \theta}\ge (1-\cos \theta) \,\frac{e^{\beta t}-e^{-\beta t}}{2},\quad \text{for $\theta\in[-\pi,\pi]$}.
\]
Thus, we have, for $t\ge 1$, that
\[
e^{\beta t \cos \theta }-e^{\beta t}\le (\cos \theta -1)\, \frac{e^{\beta t}-e^{-\beta t}}{2}
\le e^{\beta t}\, \frac{1-e^{-2 \beta}}{2}\, (\cos\theta -1)
\le -\frac{1-e^{-2 \beta}}{\pi^2}\, e^{\beta t} \,\theta^2.
\]
In the second
 inequality we have used that $t\ge 1$; the inequality $1-\cos x\ge 2x^2/\pi^2$, valid for $x\in[-\pi,\pi]$, yields the last step.
\end{proof}

%

\begin{theo}\label{teor:general for quasiexponentials} Let $g(z)=\sum_{n=0}^{\infty} b_n z^n$ satisfy
\begin{equation}\label{eq:quasiexponential}
B \,\frac{\beta^n}{n!} \le b_n\le L\,\frac{\lambda^n}{n!} \quad \mbox{for $n \ge 1$}\,,
\end{equation}
for some constants $B, L>0$, and $\beta, \lambda\ge 0$ such that $2\lambda <3\beta $.

Then $g$ satisfies the requirements of Theorem {\rm\ref{teor:cut for e^g in Hayman class bis}} and therefore $f=e^g$ is in the Hayman class.
\end{theo}
\begin{proof}
Condition  \eqref{eq:quasiexponential} gives that the power series $g$ has radius of convergence $R=+\infty$. Besides, condition \eqref{eq:quasiexponential} and Proposition \ref{prop:binomicos factoriales} show that the variance $\sigma^2(t)$ satisfies that $t^2 e^{\beta t}=O(\sigma^2(t))$ and that $\sigma^2(t)=O( t^2 e^{\lambda t})$, and that $\omega_g$ satisfies that $\omega_g(t)=O(t^3 e^{\lambda t})$,  as $t \to \infty$.

\smallskip

The variance condition \eqref{eq:variance condition for e^g} is thus obviously satisfied. The major arc condition~\eqref{eq:arco mayor para f=e^g} holds with the cut function $h(t)=e^{-\alpha t}$, where $\alpha\in(\lambda/3,\beta/2)$, since $3\alpha>\lambda$.

\smallskip

Also, for  $t >0$ and $|\theta|\le \pi$, we have that
\[
\begin{aligned}
\Re g(te^{\imath \theta})-g(t)&=\sum_{n=1}^\infty b_n\, t^n (\cos n \theta -1)\le B \sum_{n=1}^\infty \frac{(\beta t)^n}{n!} (\cos n \theta -1)\\&=B\big(\Re e^{\beta t e^{\imath \theta}}-e^{\beta t}\big)\le B  ( |e^{\beta t e^{\imath \theta}} |-e^{\beta t} )=B  (e^{\beta t \cos \theta }-e^{\beta t} )\,,
\end{aligned}
\]
and, consequently, condition \eqref{eq:major arc for quasiexponentials} of Lemma \ref{lema:major arc for quasiexponentials} is satisfied.

Condition \eqref{eq:arco menor para f=e^g bis uno uno} is satisfied with the choice $V(t)=Ce^{\beta t}$, $t>1$, and \eqref{eq:arco menor para f=e^g bis dos dos}  holds because $2\alpha<\beta$. As mentioned above, these two conditions imply  the conditions \eqref{eq:arco menor para f=e^g bis uno} and \eqref{eq:arco menor para f=e^g bis dos} of Theorem \ref{teor:cut for e^g in Hayman class bis}, and thus we conclude that $f$ is in the Hayman class.\end{proof}


\subsection{Case \texorpdfstring{$g$}{g} not entire, \texorpdfstring{$R<\infty$}{R<infinty}}

Here, we will  repeatedly appeal to Proposition  \ref{prop:binomicos}.

\begin{lem}\label{lema:major arc for quasigeometric}
If a power series $g$ with nonnegative coefficients and radius of convergence $R=1$ satisfies, for some $\beta >0$ and $B>0$, that
\begin{equation}\label{eq:major arc for quasigeometric}
\Re g(t e^{\imath \theta})-g(t)\le B \Big(\frac{1}{|1-te^{\imath \theta}|^\beta} -\frac{1}{(1-t)^\beta}\Big), \quad  \mbox{for $t\in (0,1)$ and $|\theta|\le \pi$}\,,
\end{equation}
then
\[
\sup\limits_{|\theta|\ge \omega} \Re g(te^{\imath \theta})-g(t)\le - C \, \frac{1}{(1-t)^{2+\beta}}\, \omega^2\, , \quad \mbox{for $t \in (1/2, 1)$ and $0\le\omega\le D (1-t)$}\,,\]
where $C>0$ and $D>0$  depend only on $\beta$ and $B$.

And, in particular, if we set  $V(t)= C/(1-t)^{2+\beta}$, and $U(t)=D(1-t)$, for $t\in (1/2,1)$, then  condition \eqref{eq:arco menor para f=e^g bis uno} is satisfied and  condition \eqref{eq:arco menor para f=e^g bis dos} requires that
\[\lim_{t \to \infty} \sigma(t) \exp\big(- C  h(t)^2/(1-t)^{2+\beta}\big)=0\,.\]
\end{lem}

\begin{proof}
Let $\beta>0$. For $t \in (0,1)$ and $|\theta|\le \pi$, we have that
\[
\Re g(t e^{\imath \theta})-g(t)\le \frac{B}{(1-t)^\beta} \Big(\Big(\Big|\frac{1-t}{1-te^{\imath \theta}}\Big|^2\Big)^{\beta/2}-1\Big)\,.
\]
Thus, for $t\in (1/2, 1)$ and $0\le \omega <(1-t)$ we have that
\[
\sup\limits_{|\theta|\ge \omega} \Re g(t e^{\imath \theta})-g(t)\le \frac{B}{(1-t)^\beta}
\Big(\big(1-C \frac{1}{(1-t)^2}\, \omega^2\big)^{\beta/2}-1\Big)\,,
\]
and, for $\omega <D_\beta (1-t)$, for $D_\beta>0$ appropriately small, and some  constant $C_\beta>0$, we have that
\[
\sup\limits_{|\theta|\ge \omega} \Re g(t e^{\imath \theta})-g(t)\le -\frac{B}{(1-t)^\beta}\,
C_\beta \,\frac{1}{(1-t)^2}\, \omega^2\,.
\qedhere\]
\end{proof}

\begin{theo}\label{teor:general for quasigeometric} Let $g(z)=\sum_{n=0}^{\infty} b_n z^n$ satisfy
\begin{equation}\label{eq:quasigeometric}
B \,\frac{n^\beta}{R^n} \le b_n\le L \,\frac{n^\lambda}{R^n} \, , \quad \mbox{for $n \ge 1$}\,,
\end{equation}
for some constants $B, L>0$, finite radius $R>0$ and $\beta, \lambda>-1$ such that $ 2\lambda <3\beta +1$.

Then $g$ satisfies the requirements of Theorem {\rm\ref{teor:cut for e^g in Hayman class bis}} and therefore $f=e^g$ is in the Hayman class.
\end{theo}

\begin{proof} Because of \eqref{eq:quasigeometric}, the power series $g$ has radius of convergence $R$. By considering $g(Rz)$, we may assume that $R=1$.
Appealing to Proposition  \ref{prop:binomicos}, we see that the variance function of $f=e^g$ satisfies that ${1}/{(1-t)^{3+\beta}}=O(\sigma^2(t))$ and that $\sigma^2(t)=O( {1}/{(1-t)^{3+\lambda}})$, and that $\omega_g$ satisfies that $\omega_g(t)=O( {1}/{(1-t)^{4+\lambda}} )$ as $t \uparrow 1$.

\smallskip

The variance condition \eqref{eq:variance condition for e^g} (or directly \eqref{eq:condicion de varianza en Hayman}) is  obviously satisfied.

\smallskip

We propose a cut $h(t)=(1-t)^{\alpha}$, where $\alpha \in (\lambda/3+4/3, \beta/2+3/2)$.
The major arc condition \eqref{eq:arco mayor para f=e^g} holds since $3\alpha >\lambda+4$.

For $t \in (0,1)$ and $|\theta|\le \pi$, we have that
\[\begin{aligned}
\Re g(t e^{\imath \theta})&-g(t)=\sum_{n=1}^{\infty} b_n \,t^n (\cos n\theta-1)\le B\sum_{n=1}^{\infty} n^\beta\, t^n (\cos n\theta-1)\\
&\le B_\beta \sum_{n=1}^{\infty} \frac{\Gamma(n+\beta+1)}{\Gamma(\beta+1) n!}\, t^n (\cos n\theta-1)
=B_\beta \Big(\Re \frac{1}{(1-te^{\imath \theta})^{\beta+1}}-\frac{1}{(1-t)^{\beta+1}}\Big)\\&\le B_\beta \Big(\Big|\frac{1}{1-te^{\imath \theta}}\Big|^{\beta+1}-\Big(\frac{1}{1-t}\Big)^{\beta+1}\Big)\,,
\end{aligned}
\]
where we have appealed to the comparison \eqref{eq:comparison coeffs of (1-z)beta} within Proposition \ref{prop:binomicos}. Consequently, condition \eqref{eq:major arc for quasigeometric} of Lemma \ref{lema:major arc for quasigeometric} does hold.

If we set  $V(t)= C/(1-t)^{\beta+3}$, for $t \in (1/2,1)$, \eqref{eq:arco menor para f=e^g bis uno} is satisfied and~\eqref{eq:arco menor para f=e^g bis dos} requires that
$\lim_{t \to \infty} \sigma(t) \exp(- C h(t)^2/(1-t)^{\beta+3})=0$, which does hold since $2 \alpha<\beta+3$.

Thus the conclusion follows from Theorem \ref{teor:cut for e^g in Hayman class bis}.
\end{proof}

\subsection{Applications of Theorems \ref{teor:general for quasiexponentials} and \ref{teor:general for quasigeometric}}

\subsubsection{Sets of labeled classes} \mbox{}

\medskip

$\bullet$ The egf $e^{e^z-1}$ of the Bell numbers, enumerating \textit{sets of sets},  is $f=e^g$ with $b_n=1/n!$, for $n \ge 1$.

The egf $e^{ze^z}$ of \textit{sets of pointed sets}, is $f=e^g$ with $b_n=n/n!$, for $n \ge 1$.

Both satisfy the  hypothesis of Theorem \ref{teor:general for quasiexponentials} and, thus, in particular they are in the Hayman class.

\medskip

$\bullet$ The power series $f(z)=e^{z/(1-z)}$, egf of sets of lists,  has $g(z)=z/(1-z)$ with $b_n=1$, for $n \ge 1$. The coefficients of the function $g$ satisfy the  hypothesis of Theorem~\ref{teor:general for quasigeometric} and, thus, in particular, $f=e^g$ is in the Hayman class.

In general, $f(z)=e^{z/(1-z)^\gamma}$ is in the Hayman class for $\gamma>0$.

But, $f(z)=e^{\ln 1/(1-z)}=1/(1-z)$, the egf of \textit{sets of cycles}  is not even Gaussian, as shown in Section \ref{section:gaussianity exponentials}.

\medskip

$\bullet$ The egf $f$ of \textit{sets of functions} is $f=e^g$ with $g(z)=\sum_{n=1}^\infty (n^n/n!) z^n$. In this case, the coefficients $b_n$ of $g$ satisfy \eqref{eq:quasigeometric} with $R=1/e$ and $\beta=\lambda=-1/2$, and therefore~$f$ is in the Hayman class.

In general, the exponential of $\sum_{n=1}^\infty (n^{n{-}\alpha}/n!) \,z^n$ is in the Hayman class if $0\le \alpha<1/2$. As we have seen   in Section \ref{section:gaussianity exponentials}, for $\alpha =1/2$ it is not strongly Gaussian, and for $\alpha>1/2$ is not even Gaussian.

\subsubsection{Sets of unlabeled classes} \mbox{}

\medskip

 $\bullet$ For  the ogf $P$ of partitions,  we have that $P=e^g$, with
\[
g(z)=\sum_{n=1}^\infty \frac{\sigma_1(n)}{n}\,  z^n\,,
\]
where $\sigma_1(n)$ denotes the sum of divisors of the positive integer $n$.
The coefficients $b_n=\sigma_1(n)/n$ satisfy, in this case,
\[
1\le b_n \le D_\varepsilon n^\varepsilon\,,\]
for each $\varepsilon >0$ and some constant $D_\varepsilon>0$. See Theorem 322 in \cite{HardyWright}. (Actually, we may bound from above with $\ln \ln n$, see Theorem 323 in~\cite{HardyWright}). Therefore the  function~$g$ satisfies the  hypothesis of Theorem \ref{teor:general for quasigeometric} and this means, in particular, that $P=e^g$ is in the Hayman class.

The same argument gives that the infinite product  $\prod_{j=1}^\infty  1/(1-z^j)^{c_j}$, where the $c_j$ are integers satisfying $1\le c_j\le c$, for some constant $c>1$, is in the Hayman class.

\medskip

$\bullet$ For the ogf $Q$ of partitions into distinct parts, given by
\[Q(z)\stackrel{(1)}{=}\prod_{j=1}^\infty (1+z^j), \quad \mbox{or alternatively, by} \quad Q(z)\stackrel{(2)}{=}\prod_{j=0}^\infty\frac{1}{1-z^{2j+1}},
\]
we have that  $Q=e^g$, where
\[g(z)\stackrel{\mbox{\tiny by $(1)$}}{=}\sum_{k,j\ge 1}\frac{(-1)^{k+1}}{k} \,z^{kj}\stackrel{\mbox{\tiny by $(2)$}}{=}\sum_{k\ge 1; j\ge 0}\frac{1}{k} \,z^{k(2j+1)}=\sum_{n=1}^\infty \frac{\sigma^{\text{odd}}_1(n)}{n} \, z^n\,,\]
where $\sigma^{\text{odd}}_1(n)$ registers the sum of the odd divisors of $n$.

The coefficients $b_n=\sigma^{\text{odd}}_1(n)/n$ satisfy
\begin{equation}\label{eq:bounds sum odd divisors}\frac{1}{n}\le b_n \le D_\varepsilon  n^\varepsilon\, ,\end{equation}
for each $\varepsilon >0$ and some constant $D_\varepsilon>0$. The inequality on the left holds simply because $1\mid n$, while the inequality on the right holds because $\sigma^{\text{odd}}_1(n)\le \sigma_1(n)$ and  Theorem 322 in \cite{HardyWright}.

These bounds, though, are not within reach of Theorem \ref{teor:general for quasigeometric}, and for the function~$Q$ we will verify directly the hypothesis of Theorem \ref{teor:cut for e^g in Hayman class bis}.

Using \eqref{eq:bounds sum odd divisors}, Proposition \ref{prop:binomicos} and formula \eqref{eq:formula for sigma in terms of g}, we have for the variance function $\sigma^2(t)$ of $f=e^g$ that
\[\dfrac{1}{(1-t)^{2}}=O(\sigma^2(t)) \quad \mbox{and} \quad \sigma^2(t)=O\Big(\dfrac{1}{(1-t)^{3+\varepsilon}}\Big)\,, \quad \mbox{ as $t \uparrow 1$}\,,\] and for the function $\omega_g(t)$ that  $\omega_g(t)=O({1}/{(1-t)^{4+\varepsilon}})$, as $t \uparrow 1$.

\medskip

The variance condition \eqref{eq:variance condition for e^g} (or directly \eqref{eq:condicion de varianza en Hayman}) is  obviously satisfied.
We propose a cut $h(t)=(1-t)^\alpha$, where $\alpha \in (4/3, 3/2)$.

\smallskip

The major arc condition \eqref{eq:arco mayor para f=e^g} holds since $3\alpha >4+\varepsilon$, for appropriate small $\varepsilon$.

\smallskip

For $t \in (0,1)$ and $|\theta|\le \pi$,  and $z=te^{\imath \theta}$, we have that
\[\Big|\frac{1+te^{\imath \theta}}{1+t}\Big|^2=1+\frac{2t (\cos \theta -1)}{(1+t)^2}\le 1+\frac{t}{2} (\cos \theta -1)\le e^{(t/2 (\cos \theta -1)}=e^{(\Re z -|z|)/2}\, ,\]
and, so,
\[\frac{|Q(z)|}{Q(|z|)}\le \exp\Big( \frac{1}{4} \Big(\Re \frac{z}{1-z} -\frac{|z|}{1-|z|}\Big)\Big)= \exp\Big( \frac{1}{4} \Big(\Re \frac{1}{1-z} -\frac{1}{1-|z|}\Big)\Big)\, ,\]
and, consequently,
\[\Re g(z) -g(|z|)=\ln \frac{|Q(z)|}{Q(|z|)} \le \frac{1}{4} \Big(\Re \frac{1}{1-z} -\frac{1}{1-|z|}\Big)\le\Big|\frac{1}{1-z}\Big| -\frac{1}{1-|z|}\cdot\]

Lemma \ref{lema:major arc for quasigeometric} gives us $V(t)=C/(1-t)^3$ and $U(t)=D(1-t)$. Now, condition \eqref{eq:arco menor para f=e^g bis dos} is satisfied since $\alpha<3/2$.

Thus, the function $g$ satisfies the conditions of Theorem \ref{teor:cut for e^g in Hayman class bis}, and, in particular, $Q=e^g$ is in the Hayman class.

\medskip

$\bullet$ For the ogf $M$ of plane partitions (see \cite{Flajolet}, p. 580),
\[
M(z)=\prod_{j=1}^\infty \frac{1}{(1-z^j)^j}\,,
\]
we have that $M=e^g$, with $g$ given  by
\[
g(z)=\sum_{n=1}^\infty \frac{\sigma_2(n)}{n} \,z^n\,,
\]
where $\sigma_2(n)$ denotes the sum of the squares of the divisors of the integer $n \ge 1$.
For each $\varepsilon>0$, there is a constant $C_\varepsilon>0$ such that
\[
n\le \frac{\sigma_2(n)}{n}\le C_\varepsilon\, n^{1+\varepsilon}\, , \quad \mbox{for each $n \ge 1$}\,.
\]
This follows since $\sigma_2(n)\le n \sigma_1(n)$ and $\sigma_1(n)\le C_\varepsilon n^{1+\varepsilon}$.
Thus we see that $g$ satisfies the conditions of Theorem \ref{teor:general for quasigeometric}, and, in particular, we obtain that $M=e^g$ is in the Hayman class.

Likewise, and more generally, we see that for integer $c\ge 0$, the ogf of colored partitions,
\[
\prod_{j=1}^\infty \frac{1}{(1-z^j)^{j^c}}\,,
\]
where each part $j$ appears in $j^c$ different colors, is in the Hayman class. Observe that $n^c\le \sigma_{c+1}(n)/n\le C_\varepsilon \,n^{c+\varepsilon}$, for $n \ge 1$.

\begin{question}
For the  number $p_s(n)$ of partitions into squares (i.e., partitions whose parts are whole squares) we have that
\[
S(z)\triangleq\prod_{j=1}^\infty \frac{1}{1-z^{j^2}}=\sum_{n=0}^\infty p_s(n) \,z^n\,.
\]
Is  $S$ Gaussian or strongly Gaussian? See \cite{Gafni} and \cite{Vaughan}.
Theorem \ref{teor:general for quasigeometric} is not applicable: the corresponding $g$ is
\[
g(z)=\sum_{m=1}^\infty \frac{1}{m}\Big(\sum_{j^2\mid m} j^2\Big) \,z^m\,,
\]
and $\sum_{j^2\mid m} j^2=1$ for any $m$ which is a product of distinct primes.

\smallskip

Same question for
\[
\prod_{j=1}^\infty \frac{1}{1-z^{2^j}}\,,
\]
whose coefficients count the number of partitions whose parts are powers of 2.
See~\cite{DeBruijn-mahler}. Theorem \ref{teor:general for quasigeometric} is not applicable, either.
\end{question}

\begin{question}
The infinite product $\prod_{j=1}^\infty (1+z^j)^{j^c}$ has nonnegative coefficients. Is it Hayman?
\end{question}

\section{On asymptotic formulas of coefficients}\label{seccion:asymptotic formulas coefficients}

Once you know  that a power series $f$ is in the Hayman class, and thus that $f$ is strongly Gaussian, you may use the asymptotic formula of Hayman of Theorem~\ref{teor:hayman asymptotic formula}, or even better, the asymptotic formula of B\'{a}ez-Duarte of Theorem \ref{teor:baez-duarte}, to yield an asymptotic formula for the coefficients of $f$.

As  a token of the general approach,  consider the partition function  $P$. Hayman's formula gives us that
\begin{equation}\label{eq:Hay-partition}
p(n)\sim \frac{1}{\sqrt{2\pi}} \,\frac{P(t_n)}{t_n^n\, \sigma(t_n)}\, , \quad \mbox{as $n \to \infty$}\,,
\end{equation}
where $m(t)$ and $\sigma^2(t)$ are the mean  and variance functions of $P$ and $t_n$ is such that $m(t_n)=n$.

As such, this formula \eqref{eq:Hay-partition} is too implicit, of course, on three counts:  $t_n$ is not explicit, and hardly ever is, and, besides, the formula involves $\sigma(t_n)$ and $P(t_n)$.

We proceed by appealing to the B\'{a}ez-Duarte formula \eqref{eq:formula de hayman-baez-duarte}, and by obtaining, on the one hand, asymptotic formulas for $\tau_n$ and $\widetilde{\sigma}(\tau_n)$, and for $P(\tau_n)$ on the other hand, as follows (see the details in Section 6.3.1 of \cite{K_uno}).

\medskip

\noindent \textit{Concerning $\tau_n$ and $\sigma(\tau_n)$}. By means of Euler's summation, we may approximate the mean function  $m(t)$ by
\[
m(e^{-s})\sim \frac{\zeta(2)}{s^2}\triangleq \widetilde{m}(e^{-s})\,, \quad \mbox{as $s\downarrow 0$,}\]
and $\sigma^2(t)$ by
\[
\sigma^2(e^{-s})\sim\frac{2\zeta(2)}{s^3}\triangleq \widetilde{\sigma}^2(e^{-s})\, , \qquad \mbox{as $s \downarrow 0$}\,,
\]
while checking  that condition \eqref{eq:condition baez-duarte} is satisfied.

Now if $\tau_n$ is given by $\tau_n=e^{-s_n}$ and $s_n=\sqrt{\zeta(2)/n}$ so that $\widetilde{m}(\tau_n)$, we may appeal to the formula of B\'{a}ez-Duarte \eqref{eq:formula de hayman-baez-duarte} and write
\begin{equation}\label{eq:Hay-partition2}
p(n)\sim\frac{1}{\sqrt{2\pi}}\,\frac{P(\tau_n)}{\tau_n^n \,\widetilde{\sigma}(\tau_n)}\,, \quad \mbox{as $n \to \infty$}\,.
\end{equation}

\medskip

\noindent \textit{Concerning $P(\tau_n)$.} Finally, Euler's summation again gives that
\[
P(e^{-s})\sim\frac{1}{\sqrt{2\pi}} \,\sqrt{s} \,e^{\zeta(2)/s}\, , \quad \mbox{as $s \downarrow 0$}\,,
\]
and thus
\[
P(e^{-s_n})\sim\frac{1}{\sqrt{2\pi}} \,\sqrt{s_n} \,e^{\zeta(2)/s_n}\, , \quad \mbox{as $n \to \infty$}\,.
\]

Upon substitution in \eqref{eq:Hay-partition2}, the Hardy--Ramanujan asymptotic formula for partitions follows:
\[
p(n)\sim \frac{1}{4\sqrt{3}} \,\frac{1}{n} \, e^{2\sqrt{\zeta(2)}\sqrt{n}}\, ,\quad \mbox{as $n\to \infty$}\,.\]

In short, this is the basic procedure for obtaining asymptotic formulas for coefficients of power series in the Hayman class. For a number of detailed examples, including the partition function above, general set constructions, and a few others, we refer the reader again to \cite{K_uno}.

\medskip

\noindent\textbf{Funding information.} Research of J.\,L. Fern\'{a}ndez and P. Fern\'{a}ndez supported by Fundaci\'{o}n Akusmatika. Research of V. Maci\'{a} was partially funded by grant MTM2017-85934-C3-2-P2 of Ministerio de Econom\'{\i}a y Competitividad of Spain and European Research Council Advanced Grant 834728.

\end{document}